\documentclass[a4paper,11pt]{article}

\usepackage{amssymb,amsmath,latexsym,amsthm}

\usepackage{xcolor}  
\usepackage{url}
\usepackage{hyperref}
\usepackage{cleveref}
\usepackage{enumerate}
\usepackage[square,sort,comma,numbers]{natbib}
\usepackage{dsfont, esint}
\usepackage{diagbox}
 \usepackage{geometry}
\usepackage{graphicx}
\usepackage{colortbl}
\usepackage{bbm}
\usepackage{mathtools}
\usepackage[title]{appendix}

\usepackage{endnotes}

\renewcommand{\textbf}[1]{\begingroup\bfseries\mathversion{bold}#1\endgroup}

\usepackage{breqn}

\usepackage{fancyhdr}
\usepackage{pstricks,pst-plot,pst-node,pstricks-add}
\usepackage{auto-pst-pdf}

\newtheorem{thm}{Theorem}[section]

\newtheorem{prop}[thm]{Proposition}

\newtheorem{lemma}[thm]{Lemma}

\theoremstyle{definition}

\newcommand{\R}{\mathbb R}
\newcommand{\Q}{\mathbb Q}
\newcommand{\Z}{\mathbb Z}
\newcommand{\N}{\mathbb N}

\numberwithin{equation}{section}

\def\XXint#1#2#3{{\setbox0=\hbox{$#1{#2#3}{\int}$}
    \vcenter{\hbox{$#2#3$}}\kern-.5\wd0}}

\allowdisplaybreaks

\makeatletter
\def\blfootnote{\xdef\@thefnmark{}\@footnotetext}
\newcommand{\subjclass}[2][2020]{%
  \let\@oldtitle\@title%
  \gdef\@title{\@oldtitle\footnotetext{#1 \emph{Mathematics subject classification.} #2}}%
}
\makeatother

\date{date}
\begin{document}

\title{A remark on dimensionality reduction in discrete subgroups}
\author{Rodolfo Viera}

\newcommand{\Addresses}{{
  \bigskip
  \footnotesize
\textsc{Universidad de La Serena, Departamento de Matemáticas, Avda. Juan Cisternas 1200, La Serena, Chile.}\par\nopagebreak
\textit{E-mail address: }\texttt{rodolfo.vieraq@userena.cl}
}}

\date{\today}
\maketitle
\begin{abstract}
In this short note, we prove a version of the Johnson-Lindenstrauss flattening Lemma for point sets taking values in discrete subgroups. More precisely, given $d,\lambda_0,N_0\in\N$ and $\epsilon\in \left(0,\frac{1}{2}\right)$ suitably chosen, we show there exists a natural number $k=k(d,\epsilon)=O\left(\frac{1}{\epsilon^2}\log d\right)$, such that for every sufficiently large scaling factor $\lambda\in\N$ and any point set $\mathcal{D}\subset\frac{\lambda}{\lambda_0}\Z^d\cap B(0,\lambda N_0)$ with cardinality $d$, there exists an embedding $F:\mathcal{D}\to\frac{1}{\lambda_0}\Z^k$, with distortion at most $\left(1+\epsilon+\frac{\epsilon}{\lambda\lambda_0}\right)$.

\end{abstract}

\maketitle

\section{Introduction}
The renowned Johnson-Lindenstrauss Lemma \cite{dimred, JLMat, Ost} (JL-Lemma for short) establishes that, given a point set $\mathcal{D}=\{x_1,\ldots,x_d\}\subset\R^d$ and a positive number $\epsilon\in (0,1)$, there exists a (linear) embedding $\Phi:\R^d\to\R^k$, where $k=k(d,\epsilon)=O(\log d/\epsilon^2)$, that maps $\mathcal{D}$ into $\R^k$ with distortion at most $(1+\epsilon)$. Although Johnson and Lindenstrauss proved their lemma to tackle a problem concerning extensions of Lipschitz maps, the computer and data science communities realized the potential of this lemma in reducing the dimension of high-dimensional data while preserving its key features up to a constant multiplicative error $\sim 1$; we refer the reader to \cite{FrJL, introJL, JLpol}  and the references therein for a broader discussion.
\medskip

Albeit JL-Lemma has become a powerful tool for "flattening" high-dimensional data, represented by vectors in $\R^d$ for some $d\gg 1$, without distorting the distances too much, in the author's opinion a more realistic scenario for a computer model-space of $d$-dimensional vectors is the set $\frac{1}{\lambda_0}\Z^d\cap B_{N_0}$, where $\lambda_0, N_0$ are fixed positive integers and $B_{N_0}$ denotes the euclidean ball with radius $N_0$ centered at the origin; this is because we cannot consider vectors with arbitrarily large entries (in absolute value) or as many decimals as we want. Thus, an interesting question is, fixed some suitable error term $\epsilon\in (0,1)$, whether we can reduce the number of variables of data in $\frac{1}{\lambda_0}\Z^d\cap B_{N_0}$ by embedding them into the grid $\frac{1}{\lambda_0}\Z^k\cap B_N$ in the same spirit as JL-Lemma, for some positive integer $N$, and $k\lesssim_{\epsilon}\log d$. 
\medskip

A naive approach is to proceed as follows: after applying the JL-Lemma to a data point set $\mathcal{D}\subset\frac{1}{\lambda_0}\Z^d\cap B_{N_0}$ with $|\mathcal{D}|=d$, we obtain a point set $\mathcal{D}_{\mathsf{flat}}=\{y_1,\ldots,y_k\}\subset\R^k$, where $k$ is given as in the conclusions of JL-Lemma; then define $\widetilde{\mathcal{D}}_{\mathsf{flat}}=\{\widetilde{z_1},\ldots,\widetilde{z_d}\}$, where $\widetilde{z_i}$ is the closest point of $\frac{1}{\lambda_0}\Z^k$ from $y_i$. A priori, for any $\lambda>0$ we can only ensure that

\[
d\left(\lambda\mathcal{D}_{\mathsf{flat}},\frac{1}{\lambda_0}\Z^k\right)\leq\frac{\sqrt{k}}{\lambda_0}.
\]

Even in the best case, by uniform distribution modulo 1 (see the proof of Lemma \ref{lem: tanerinte} below), we could find a sequence $(n_l)_{l\geq 1}$ such that 

\[
(\forall l\geq 1):\qquad d\left(n_l\mathcal{D}_{\mathsf{flat}},\frac{1}{\lambda_0}\Z^k\right)<\epsilon.
\]
\medskip

 In view of the previous discussion, we aim to prove the following: given $d,\lambda_0\in\N$ and $\epsilon\in (0,1/2)$ appropriately chosen, then {\em every} sufficiently separated point set in $\frac{1}{\lambda_0}\Z^d$ can be flattened into in a $k$-dimensional sub-lattice $\frac{1}{\lambda_0}\Z^k$ with distortion $\sim_{\epsilon,\lambda_0} 1$ (i.e., with distortion factor close 
to 1 and depending on $\epsilon$ and $\lambda_0$), for some $k=k(d,\epsilon)\ll d$.

\begin{prop}[Main Proposition]\label{prop: JLdis}
    Let $d,\lambda_0,N_0\in\N$ and $\epsilon\in\left(0,\frac{1}{\lambda_0+1}\right)$ be given. There exists $c=c(\epsilon)>0$ such that the following holds for every positive integer $k\geq c\log d$: there is a scaling factor $\lambda_1=\lambda_1(\lambda_0,\epsilon,k,N_0)\in\N$ such that for every $\lambda\geq\lambda_1$ and for every point set $\mathcal{D}\subset\frac{\lambda}{\lambda_0}\Z^d\cap B_{\lambda N_0}$ with $|\mathcal{D}|=d$, there is a mapping $F:\mathcal{D}\to\frac{1}{\lambda_0}\Z^k$ such that
    \begin{equation}\label{eq: JLdis}
        (\forall x,y\in\mathcal{D}):\qquad \left(1-\epsilon-\frac{\epsilon}{\lambda\lambda_0}\right)\|x-y\|\leq\|F(x)-F(y)\|\leq\left(1+\epsilon+\frac{\epsilon}{\lambda\lambda_0}\right)\|x-y\|.
    \end{equation}
\end{prop}

This version of JL-Lemma has the advantage that, after rescaling a data set by $\lambda$, we can reduce its dimension with small distortion while keeping the number of decimals and the magnitude of the flattened data bounded.  
\medskip

Besides the JL-Lemma itself, the following remarkable Theorem due to Tamar Ziegler \cite[Theorem 1.3]{Nil} (see also \cite{appol} for an elementary proof in the case of even dimensions) plays a crucial role in the proof of the \nameref{prop: JLdis}.

\begin{thm}[Ziegler's Theorem]\label{thm: nil}
    Let $\mathcal{D}:=\{x_1,\ldots,x_n\}$ be a set of $n$ vectors in $\R^d$ and $\epsilon>0$ be given. Then there is $l_0=l_0(\epsilon,\mathcal{D})>0$ such that for any $l\geq l_0$, there exists a rotation $\rho=\rho(l)\in SO(d)$ satisfying that
    \[
    (\forall i=1,\ldots,n):\qquad d(\rho(l\cdot x_i),\Z^d)\leq\epsilon.
    \]
\end{thm}
\section{Proof of Main Proposition}
The following result is a key ingredient to determine the dependence of the parameter $\lambda$.
\begin{lemma}
    Let  $t\in (0,1)\setminus\Q$, $N\in \N$, and $\epsilon>0$ be fixed. Then there exist $\lambda_1=\lambda_1(t,\epsilon, N)\in\N$ such that for every $\lambda\geq\lambda_1$ and $\mathcal{D}\subset t\Z^k\cap B_N$, there must exists a rotation $\rho=\rho(\lambda)\in SO(k)$ such that
    \[
    d(\rho(\lambda\mathcal{D}),\Z^k)<\epsilon.
    \]
\end{lemma}

\begin{proof}\label{lem: tanerinte}
   Since $t\in (0,1)\setminus\Q$, the sequence $(n t)_{n\geq 1}$ is uniformly distributed modulo 1 (see for instance \cite{unifseq}); in particular, given $\epsilon>0$, there exist $n_1=n_1(\epsilon,N)\in\N$ and $p\in\Z$ such that
    \[
    \left|n_1t-p\right|<\frac{\epsilon}{N}\leq\epsilon.
    \]

    Thus, for every $l=1,\ldots,k$ and $q=(q_1,\ldots,q_k)\in\Z^k\cap B_N$, we obtain:
    \begin{equation}\label{eq: tamarinte}
    \left|n_1tq_l-pq_l\right|<\frac{\epsilon|q_l|}{N}\leq\epsilon
     \end{equation}
    Hence from \eqref{eq: tamarinte}, for every subset $\mathcal{D}$ of $t\Z^k\cap B_N$ there holds that $d(n_1\mathcal{D},\Z^k)<\epsilon$. The rest of the proof follows the very same lines as in \cite{appol}.
\end{proof}

\begin{proof}[Proof of the \nameref{prop: JLdis}]
    Let $\mathcal{D}_{\mathsf{flat}}:=\Phi(\mathcal{D})$, where $\Phi:\R^d\to\R^k$ is the linear embedding given by the Johnson-Lindenstrauss Lemma, and write $y_i:=\Phi(x_i)$ where $\mathcal{D}=\{x_1,\ldots,x_d\}$; the proof is quite direct if $k$ is a perfect square, and so we assume that $\sqrt{k}$ is an irrational number. Since a translation by a vector is an isometry, we can assume that the origin of $\R^k$ is the circumcenter of $\mathcal{D}_{\mathsf{flat}}$. Moreover, by \cite[Theorem 3.1]{JLMat} (see also \cite[Theorem 1.35]{Ost}), we can consider $\mathcal{D}_{\mathsf{flat}}$ as a subset of $\frac{1}{\lambda_0\sqrt{k}}\Z^k$, since $\Phi$ takes the form:
\[
\Phi(x)=\frac{1}{\sqrt{k}}Rx^{\mathsf{T}},
\]

where $R$ is a $d\times k$ (random matrix) with entries taking values in $\{0,1\}$. In particular, by linearity, for every $t>0$ we have that
\[
t\mathcal{D}_{\mathsf{flat}}=\Phi(t\mathcal{D}), 
\]
and in consequence, we get that $\Phi\circ\mathsf{dil}_t:\mathcal{D}\to t\mathcal{D}_{\mathsf{flat}}\subset\R^k$ is an $(1+\epsilon)$-embedding, where $\mathsf{dil}_t$ stands for the dilation by $t$.
\medskip

By \nameref{thm: nil} and Lemma \ref{lem: tanerinte}, there exists a scaling factor $\lambda_1=\lambda_1(\epsilon,k, N_0)\in\N$ such that for every $\lambda\geq\lambda_1$ there is a rotation $\rho=\rho(\lambda)\in SO(k)$, such that 
\[
d\left(\rho(\lambda\mathcal{D}_{\mathsf{flat}}),\frac{1}{\lambda_0}\Z^k\right)\leq\frac{\epsilon}{\lambda_0}\leq\frac{\epsilon}{2},
\]
and thus for each $i\in\{1,\ldots,d\}$, there exists $z_i\in\Z^k$ such that
\begin{equation}\label{eq: aproxrot1}
\left\|\rho(\lambda y_i)-\frac{1}{\lambda_0}z_i\right\|<\frac{\epsilon}{2}.
\end{equation}

Firstly we shall prove that for $i\neq j$, there holds that $z_i\neq z_j$; indeed, by \eqref{eq: aproxrot1}, the linearity of $\Phi$, the triangle inequality, the definition of $\mathcal{D}_{\mathsf{flat}}$, and since the $x_i$'s belong to $\frac{1}{\lambda_0}\Z^d$, we have that 
\[
\frac{(1-\epsilon)}{\lambda_0}\leq (1-\epsilon)\|x_i-x_j\|\leq\|y_i-y_j\|\leq\frac{\epsilon}{\lambda}+\frac{1}{\lambda\lambda_0}\|z_i-z_j\|,
\]
and thus, by the choice of $\epsilon$, we have that
\[
\frac{1}{\lambda\lambda_0}\|z_i-z_j\|\geq \frac{1}{\lambda_0}-\epsilon\left(1+\frac{1}{\lambda_0}\right)>0.
\]

Now we claim that the mapping $F:\lambda x_i\mapsto\frac{1}{\lambda_0}z_i$ verifies \eqref{eq: JLdis}. By triangle inequality, the fact the JL-embedding $\Phi$ has distortion at most $(1+\epsilon)$, and that the $x_i$'s belong to $\frac{1}{\lambda_0}\Z^d$, and since $\rho\in SO(k)$ is an isometry, we get
\begin{equation}\label{eq: JLdis1}
    \frac{1}{\lambda_0}\|z_i-z_j\|\leq \epsilon + \|\rho(\lambda y_i)-\rho(\lambda y_j)\|\leq\epsilon+\lambda(1+\epsilon)\|x_i-x_j\|\leq \left(1+\epsilon+\frac{\epsilon}{\lambda\lambda_0}\right)\|\lambda x_i-\lambda x_j\|.
\end{equation}

Analogously, we have that 
\begin{equation}\label{eq: JLdis2}
    \left(1-\epsilon-\frac{\epsilon}{\lambda\lambda_0}\right)\|x_i-x_j\|\leq\frac{1}{\lambda_0}\|z_i-z_j\|.
\end{equation}

Therefore, \eqref{eq: JLdis} follows from \eqref{eq: JLdis1} and \eqref{eq: JLdis2}. This finishes the proof.

\end{proof}

\Addresses
\end{document}